\documentclass[11pt]{article}       
\usepackage{graphicx}
\usepackage{natbib}
\usepackage{bbm}
\usepackage{fullpage}
\usepackage{graphicx}
\usepackage{setspace}
\usepackage{xcolor}

\usepackage{amsmath, amsthm, amssymb}
\usepackage{tikz}
\usepackage{subfig}
\usepackage{multirow}
\usepackage{float}
\usepackage{mathtools}
\RequirePackage[colorlinks,citecolor=blue,urlcolor=blue]{hyperref}
\DeclareMathOperator*{\argmin}{argmin}

\DeclareMathOperator{\sign}{sign}
\DeclareMathOperator{\GCV}{GCV}

\DeclareMathOperator{\Tr}{Tr}

\definecolor{DR}{rgb}{0.8509804, 0.3176471, 0.3058824} 
\definecolor{TS}{rgb}{0.1647059,  0.1686275, 0.1764706}
\definecolor{BB}{rgb}{0.1764706, 0.6588235, 0.8470588}

\definecolor{c1}{rgb}{0,  1, 0}
\definecolor{c2}{rgb}{0,  0.39, 0}
\definecolor{c3}{rgb}{0.4980392, 0.4980392, 0.4980392}
\definecolor{c4}{rgb}{0.7490196, 0.2470588, 0.7490196}
\definecolor{c5}{rgb}{1, 0, 1}

\DeclareRobustCommand\full  {\tikz[baseline=-0.6ex]\draw[blue, thick] (0,0)--(0.57,0);} 
\DeclareRobustCommand\fulll  {\tikz[baseline=-0.6ex]\draw[red, thick] (0,0)--(0.57,0);} 
\DeclareRobustCommand\dotted{\tikz[baseline=-0.6ex]\draw[c2,thick,dotted] (0,0)--(0.57,0);} 
\DeclareRobustCommand\denselydashed{\tikz[baseline=-0.6ex]\draw[c5, thick, dash pattern={on 7pt off 1.5pt}] (0,0)--(0.57,0);} 
\DeclareRobustCommand\dotdash {\tikz[baseline=-0.6ex]\draw[c4,thick,dash dot] (0,0)--(0.62,0);} 
\DeclareRobustCommand\dashed {\tikz[baseline=-0.6ex]\draw[red,thick,dashed] (0,0)--(0.57,0);} 
\DeclareRobustCommand\dashedd {\tikz[baseline=-0.6ex]\draw[green,thick,dashed] (0,0)--(0.57,0);}

\newtheorem{cor}{Corollary}
\newtheorem{theorem}{Theorem}

\begin{document}

\title{Asymptotics for M-type smoothing splines with non-smooth objective functions
}


\author{Ioannis Kalogridis}


\maketitle

\begin{abstract}
M-type smoothing splines are a broad class of spline estimators that include the popular least-squares smoothing spline but also spline estimators that are less susceptible to outlying observations and model-misspecification. However, available asymptotic theory only covers smoothing spline estimators based on smooth objective functions and consequently leaves out frequently used resistant estimators such as quantile and Huber-type smoothing splines. We provide a general treatment in this paper and, assuming only the convexity of the objective function, show that the least-squares (super-)convergence rates can be extended to M-type estimators whose asymptotic properties have not been hitherto described. We further show that auxiliary scale estimates may be handled under significantly weaker assumptions than those found in the literature and we establish optimal rates of convergence for the derivatives, which have not been obtained outside the least-squares framework. A simulation study and a real-data example illustrate the competitive performance of non-smooth M-type splines in relation to the least-squares spline on regular data and their superior performance on data that contain anomalies.
\end{abstract}

\section{Introduction}
\label{sec:SMSP1}
Based on data $(t_{1}, Y_1), \ldots, (t_{n}, Y_n)$ with non-random $t_{i} \in [0,1]$, consider the classical nonparametric regression model
\begin{equation}
\label{eq:SMSP1}
Y_i = f_o(t_{i}) + \epsilon_i, \quad (i=1, \ldots, n),
\end{equation}
where $f_o$ is a sufficiently smooth function that we would like to estimate and the $\epsilon_i, \ i=1,\ldots, n$, are independent and identically distributed error terms, commonly assumed to have zero mean and finite variance $\sigma^2$.

A popular estimation method involves restricting $f_o$ to the Hilbert-Sobolev space of smooth functions denoted by $\mathcal{W}^{m,2}([0,1])$ and defined as
\begin{align*}
\mathcal{W}^{m, 2}([0,1]) = \{ f:[0,1] \to \mathbb{R}, f\ &\text{has $m-1$ absolutely continuous derivatives}  \\  & f^{(1)}, \ldots, f^{(m-1)}\ \text{and} \int_0^1 |f^{(m)}(x)|^2 dx< \infty  \},
\end{align*}
and finding $\widehat{f}_n \in \mathcal{W}^{m,2}([0,1])$ which minimizes
\begin{equation}
\label{eq:SMSP2}
\frac{1}{n} \sum_{i = 1}^n | Y_{i} - f(t_{i}) |^2 + \lambda \int_0^{1} |f^{(m)}(x)|^2 dx,
\end{equation}
for some $\lambda>0$ that governs the trade-off between smoothness and goodness-of-fit. The problem is well-defined for $n \geq m$ and its solution is a $2m$th order natural spline with knots at $t_{1}, \ldots, t_{n}$. The subsequent least-squares smoothing spline can be computed very efficiently with the Kimeldorf-Wahba representer theorem and can be shown to attain the optimal rates of convergence, under the usual Gauss-Markov conditions on the error term. The interested reader is referred to \citep{Wahba:1990}, \citep{Green:1994} and \citep{Eubank:1999} for detailed theoretical developments and illustrative examples.

The focus of this paper is theoretical, as we aim show that, under weak assumptions, $f_o$ may be optimally estimated even if the error lacks finite moments. The estimator in consideration is the M-type smoothing spline introduced by \cite{Huber:1979} and defined as a solution of 
\begin{equation}
\label{eq:SMSP3}
\inf_{f \in \mathcal{W}^{m, 2}([0,1])}\left[ \frac{1}{n} \sum_{i = 1}^n \rho\left( Y_{i} - f(t_{i}) \right) + \lambda \int_0^{1} |f^{(m)}(x)|^2 dx \right],
\end{equation}
for some convex nonnegative function $\rho$ that is symmetric about zero and satisfies $\rho(0)=0$. Clearly, the least-squares smoothing spline fulfils these conditions but the benefit of the above formulation is that it allows for more general loss functions that reduce the effect of large residuals and make for resistant estimation of $f_o$. Popular examples include the resistant absolute value loss $|x|$ and Huber's function
\begin{equation*}
\rho_k(x) = \begin{cases} x^2/2, & |x| \leq k \\ k\left(|x|-k/2\right), & |x| >k
\end{cases}\\
\end{equation*}
where the tuning parameter $k$ controls the blending of square and absolute losses. It is shown below that under general conditions a solution to \eqref{eq:SMSP3} in $\mathcal{W}^{m,2}([0,1])$ exists, although, as \cite{Cox:1983} and \cite{Eg:2009} remark, it may not be unique unless $\rho$ is strictly convex. Similarly to the least-squares setting, if $n \geq m$ then this minimizer must be a $2m$th order natural spline with knots at the design points. 

Traditionally in robust regression the losses are standardized with an equivariant scale estimator $\widehat{\sigma}$ in order to achieve scale equivariance of the regression estimator. The M-type smoothing spline estimator defined in \eqref{eq:SMSP3} cannot, in general, be scale equivariant even after this standardization on account of the penalty term. Nevertheless, we may define an "approximately" equivariant estimator as the solution of 
\begin{equation}
\label{eq:SMSP4}
\inf_{f \in \mathcal{W}^{m, 2}([0,1])}\left[ \frac{1}{n} \sum_{i = 1}^n \rho\left( \frac{Y_{i} - f(t_{i})}{\widehat{\sigma}} \right) + \lambda \int_0^{1} |f^{(m)}(x)|^2 dx \right].
\end{equation}
Originally, \citet{Huber:1979} proposed simultaneous scale estimation but nowadays it is recognized that preliminary scale estimates tend to perform better, see, e.g., \citep{Maronna:2006}. Such scale estimates may be obtained either from preliminary model fitting, that is, from fitting a robust regression estimator to the data and computing a robust scale from its residuals, or from robust scale estimates involving linear combinations of the $Y_{i}$s, \citep{Cunningham:1991, Ghement:2008}. Inclusion of $\widehat{\sigma}$ adds a new theoretical layer to the smoothing spline problem and its asymptotic properties influence those of the smoothing spline estimator.

In stark contrast to least-square smoothing splines, where numerous theoretical results have been obtained ranging from equivalent kernels to convergence rates of the derivatives, only a few works have delved into the theory of general M-type smoothing splines. \citet{Cox:1983} and \citet{Oh:2007} obtained an asymptotic linearization of M-type estimators with $\rho$ functions in $\mathcal{C}^{3}(\mathbb{R})$ and used this to show that smooth M-type estimators attain the least-squares convergence rates. \citet{Cunningham:1991} complemented the work of the first author by showing that for the special case of $m=2$ the optimal rates of convergence are retained if one uses a root-n preliminary scale provided that the error term possesses a first moment. \citet{van de Geer:2002} was able to reduce the smoothness requirements to a Lipschitz condition on $\rho$, but her work does not address either the case of auxiliary scale estimates or estimation of derivatives. Finally, \citet{Eg:2009} study the least-absolute deviations smoothing spline in detail, but only under the assumption that it is asymptotically contained in a ball around $f_o$. 

We provide a unified treatment of the M-type smoothing spline problem in this paper, including preliminary scale estimation and estimation of derivatives. Our main assumptions center around a convex loss function, two mild regularity conditions on the errors $\epsilon_i$ that have been widely used in the non-penalized case, and approximate uniformity of the design points. For well-chosen loss functions, these conditions do not require the existence of any moments of the error allowing for very
heavy-tailed distributions, under which the least-squares estimator may fail to be consistent. Furthermore, we show that these conditions barely change with the inclusion of auxiliary scale estimates constructed either from pseudo-residuals or preliminary regression estimates. Our treatment is potentially of independent mathematical interest, as it relies on the theory of reproducing kernel Hilbert spaces instead of the nowadays commonly used theory of empirical processes \citep{van de Geer:2000, van de Geer:2002}).

\section{Main results: existence of solutions and rates of convergence}
\label{sec:SMSP2}

We begin by introducing some useful notation. We denote the standard $\mathcal{L}^2([0,1])$ inner product by $\langle \cdot, \cdot \rangle_2$. Its associated norm $||\cdot||_2$ is given by $\langle f, f\rangle_2^{1/2}$ for any $f \in \mathcal{L}^2([0,1])$. Throughout, we endow $\mathcal{W}^{m,2}([0,1])$ with the inner product given by
\begin{equation}
\label{eq:SMSP5}
\langle f, g \rangle _{m, \lambda} = \langle f, g \rangle_{2} + \lambda \langle f^{(m)}, g^{(m)} \rangle_{2},
\end{equation}
for any $f, g \in \mathcal{W}^{m,2}([0,1])$. The associated norm is denoted by $||\cdot||_{m, \lambda}$. Norms depending on the smoothing parameter have also been used by \citet{Silverman:1996} and \citet{Eg:2009}, for example. Denoting the M-type smoothing spline by $\widehat{f}_n$, we shall see that an advantage of this norm is that establishing rates of convergence with respect to $||\widehat{f}_n-f_o||_{m,\lambda}$ will semi-automatically  yield rates of convergence for the derivatives with respect to $||\cdot||_{2}$. 

By an extension of the Sobolev embedding theorem, \citet{Eg:2009} showed that for all $x\in[0,1]$, all $f \in \mathcal{W}^{m,2}([0,1])$ and all $\lambda \in (0,1)$, there exists a constant $c_m$, depending only on $m$, such that
\begin{equation}
\label{eq:SMSP6}
|f(x)| \leq \frac{c_m}{\lambda^{1/4m}} ||f||_{m,\lambda}.
\end{equation}
This result implies that point-evaluation is a continuous linear functional with the inner product \eqref{eq:SMSP5}. It follows that $\mathcal{W}^{m,2}([0,1])$ is a reproducing kernel Hilbert space and  there exists a symmetric function $\mathcal{R}_{m,\lambda}(x, y)$, the reproducing kernel, such that $\mathcal{R}_{m, \lambda}(x, \cdot) \in \mathcal{W}^{m, 2}([0,1])$ for every $x \in[0,1]$ and for every $f \in \mathcal{W}^{m,2}([0,1])$,
\begin{equation}
f(x) = \langle f, \mathcal{R}_{m,\lambda}(x, \cdot) \rangle_{m,  \lambda}.
\end{equation}
Consequently, by \eqref{eq:SMSP6},
\begin{equation}
\label{eq:SMSP8}
\sup_{x \in [0,1]} ||\mathcal{R}_{m, \lambda}(x, \cdot)||_{m, \lambda} \leq \frac{c_m}{\lambda^{1/4m}},
\end{equation}
with the same constant $c_m$. The above bounds will play a key role in the establishment of our results.

We first deal with the existence of the M-type smoothing spline and show that the problem \eqref{eq:SMSP3} is well-defined, in the sense that it possesses at least one solution in $\mathcal{W}^{m, 2}([0,1])$. The theorem requires only a weak form of continuity of $\rho$ and may therefore be useful in other settings as well.
\begin{theorem}
\label{Thm:SMSP1}
If $\rho(x)$ is a lower semicontinuous, nondecreasing unbounded function of $|x|$ and $n \geq m$, the minimization problem 
\begin{equation*}
\inf_{f \in \mathcal{W}^{m, 2}([0,1])}\left[ \frac{1}{n} \sum_{i = 1}^n \rho\left( Y_{i} - f(t_{i}) \right) + \lambda ||f^{(m)}||_2^2 \right],
\end{equation*}
has a solution in $\mathcal{W}^{m,2}([0,1])$.
\end{theorem}
Arguing in a standard way now shows that a minimizer may be found in the $n$-dimensional space of  natural splines of order $2m$ with knots at $t_{1}, \ldots, t_{n}$. Commonly used convex $\rho$-functions, such as Huber's function, can be easily shown to satisfy the condition of Theorem \ref{Thm:SMSP1}. However, Theorem~\ref{Thm:SMSP1} is also applicable for certain non-convex $\rho$-functions, such as $\rho(x) = \log(1+x^2)$. As mentioned previously, existence can be strengthened to uniqueness if one uses a strictly convex $\rho$-function, such as the logistic $\rho$-function. See Proposition 2.1 of \citet{Cox:1983}.

We may now treat the asymptotics of M-type estimators with scale either known or, more realistically, not needed. The required regularity conditions on $\rho$, the error term and the design points $t_{1}, \ldots, t_{n}$ are as follows.

\begin{itemize}
\item[(A1)] The loss function $\rho(x)$ is absolutely continuous and convex with $\psi(x)$ any choice of its subgradient. 
\item[(A2)] There exist finite constants $\kappa$ and $M_1$ such that for all $x \in \mathbb{R}$ and $|y| < \kappa$,
\begin{align*}
|\psi(x+y)-\psi(x)| \leq M_1.
\end{align*}
\item[(A3)] There exists a finite constant $M_2$ such that
\begin{align*}
\sup_{|t| \leq h} \mathbb{E}\{ |\psi(\epsilon_1+t) - \psi(\epsilon_1)|^2 \} \leq M_2 |h|,
\end{align*}
as $h \to 0$.
\item[(A4)] $\mathbb{E}\{|\psi(\epsilon_1)|^2\} \leq \tau^2 < \infty$, $\mathbb{E}\{\psi(\epsilon_1)\} = 0$ and there exists a constant $\xi>0$ such that
\begin{equation*}
\mathbb{E}\left\{\psi(\epsilon_1+t) \right\} = \xi t + o(t), \quad \text{as} \quad t \to 0.
\end{equation*}
\item[(A5)] The family of the design points $\{t_{i}\}_{i=1}^n$ is quasi-uniform in the sense of \citet{Eg:2009}, that is, there exists a constant $M_3$ such that, for all $n \geq 2$ and all $f \in \mathcal{W}^{1, 1}([0,1])$,
\begin{equation*}
\left| \frac{1}{n} \sum_{i=1}^n f(t_{i}) - \int_{0}^1 f(t) dt \right| \leq \frac{M_3}{n} \int_{0}^{1}|f^{\prime}(t)| dt
\end{equation*}
\end{itemize}

Condition (A1) is standard in the asymptotic theory of M-estimators for unpenalized linear models, see, for example, \citep{Yohai:1979, Bai:1994}. Condition (A2) requires that $\psi$ has uniformly bounded local increments and is the only condition that is imposed directly on $\psi$. This needs to be contrasted with the restrictive smoothness conditions of \citep{Cox:1983}, \citep{Cunningham:1991} and \citep{Oh:2007}, all of whom assumed a twice-differentiable $\psi$-function with bounded second derivative. In the same spirit, conditions (A3)--(A4) trade differentiability of $\psi$ with some regularity of the distribution of the error term. They are very mild. Condition (A3) is a mean-square continuity condition on $m(t) := \psi(\epsilon_1+t)$ and holds quite generally. For example, \citet{Bai:1994} demonstrated that in some cases it may be possible to have the much tighter bound $M_2 |h|^2$ on the right-hand side of the inequality, even for $\psi$-functions that have jumps. See also \citep{Welsh:1989} for this point. Condition (A4) ensures the Fisher-consistency of the estimates and has been widely used for many types of M-estimators, see \citep{Huber:2009, Maronna:2006} for relevant discussions. A number of interesting estimators can be covered by condition (A4), as we now show.

Example 1 (LAD and quantile regression). Consider M-estimation with $\rho(x)=|x|$. Then, provided that $\epsilon_1$ has a distribution function $F$ that is symmetric about zero and a positive density $f$ on an interval about zero,
\begin{equation*}
\mathbb{E}\{ \sign(\epsilon_1+t) \} = 2f(0)t + o(t), \quad \text{as} \quad t \to 0.
\end{equation*}
cfr. \citep{Pollard:1991}. This generalizes easily to M-estimation with $\rho_{\alpha}(x) = |x| + (2\alpha-1)x$, provided that in this case one views the regression function as the $\alpha$-quantile function, that is, $\Pr(Y_i \leq f_o(t_i)) = \alpha$.

Example 2 (Huber). For $k>0$, $\psi_k(x) = \max\{-k, \min\{x,k\}\}$ and we may assume that $F$ is absolutely continuous and symmetric about zero so that
\begin{equation*}
\mathbb{E}\{ \psi_k(\epsilon_1+t) \} = \left( 2 F(k)-1\right)t + o(t), \quad \text{as} \quad t \to 0.
\end{equation*}
The term $2 F(k)-1$ is strictly positive for all $k>0$ under these assumptions.

Example 3 ($\mathcal{L}^{p}$ regression estimates with $ 1< p < 2$). Clearly, $\psi_p(x) = p |x|^{p-1} \sign(x)$ and if we assume that $F$ is symmetric about zero, $\mathbb{E}\{| \epsilon_1|^{p-1}\} < \infty$ and $\mathbb{E}\{|\epsilon|^{p-2} \}<\infty$ then
\begin{equation*}
\mathbb{E}\{\psi_p(\epsilon_1+t) \} = p(p-1)\mathbb{E}\{|\epsilon|^{p-2} \} t + o(t), \quad \text{as} \quad t \to 0,
\end{equation*}
see \citep{Arcones:2000}. The latter expectation is finite, if, e.g., $F$ possesses a Lebesgue density $f$ that is bounded at an interval about zero.

Example 4 (Expectile regression). As an alternative to the check loss, consider the expectile loss $\rho_{\alpha}(x) = x^2/2| \alpha - \mathcal{I}(x\leq 0)|$ with $\alpha \in (0,1)$, such that $\psi_{\alpha}(x) = (1-\alpha)x \mathcal{I}(x \leq 0) + \alpha x\, \mathcal{I}(x>0)$. Assuming that there exists an interval about zero on which $F$ has no atoms, we have
\begin{equation*}
\mathbb{E}\{\psi_{\alpha}(\epsilon_{1}+t) \} =  \{\alpha + (1- 2 \alpha) F(0)\}t + o(t), \quad \text{as} \quad t \to 0.
\end{equation*}
Therefore, assumption (A4) holds with $\xi = \alpha + (1- 2 \alpha) F(0)$ and this is bounded away from zero and infinity for all $\alpha \in (0,1)$.

Example 5 (Smooth $\rho$-functions).  All monotone everywhere differentiable score functions $\psi$ with bounded second derivative $\psi^{\prime \prime}(x)$, resulting, for example, from $\rho(x) = \log(\cosh(x))$, satisfy the second part of assumption (A4) provided that $\mathbb{E}\{\psi^{\prime}(\epsilon_{1})\}>0$. This is a classical Fisher-consistency for M-estimators based on differentiable score functions, see \citep{Maronna:2006}.

It should be noted that although for convenience (A3) and (A4) are stated here with independent and identically distributed errors in mind, these conditions can be weakened to include solely independent errors. In that case we would require some shared regularity of the errors and restate (A3) as 
\begin{align*}
\sup_{n} \max_{1 \leq i \leq n} \sup_{|t| \leq h} \mathbb{E}\{ |\psi(\epsilon_i+t) - \psi(\epsilon_i)|^2 \} \leq M_2 |h|,
\end{align*}
for some finite $M_2$, as $h \to 0$. Similarly, to extend (A4) to independent errors we would require $\sup_{n} \max_{1 \leq i \leq n} \mathbb{E}\{|\psi(\epsilon_i)|^2 \} \leq \tau^2$, $\mathbb{E}\{\psi(\epsilon_i) \}=0$, $\mathbb{E}\{\psi(\epsilon_i + t) \} = \xi_i t + o(t)$, as $ t \to 0$ for $i=1, \ldots, n$, with
\begin{align*}
0< \inf_{n} \min_{1 \leq i \leq n} \xi_i \leq \sup_n \max_{1 \leq i \leq n} \xi_i < \infty.
\end{align*}

Finally, condition (A5) ensures that the responses are observed at a sufficiently regular grid. It is a quite weak assumption that can be shown to hold for all frequently employed designs such as $t_{i} = i/n$ or $t_{i} = 2i/(2n+1)$. Call $G_n$ the distribution function that jumps $n^{-1}$ at each $t_{i}$. An integration by parts argument shows 
\begin{equation*}
\frac{1}{n}\sum_{i=1}^n f(t_{i}) - \int_{0}^1 f(x) dx = \int_{0}^1 \{x- G_n(x) \} f^{\prime}(x) dx,
\end{equation*}
and (A5) is satisfied, if, for example,  $G_n$ approximates well the uniform distribution function in the Kolmogorov metric.

It should be noted that none of the loss functions in examples 1--4 are covered by the theory of \citep{Cox:1983} and \citep{Oh:2007}, while the Lipschitz condition employed by \citet{van de Geer:2002} leaves out $\mathcal{L}^p$ and expectile smoothing spline estimators, among others. Our first asymptotic result is Theorem \ref{Thm:SMSP2}, which establishes the optimality of general M-type smoothing splines, provided that the smoothing parameter $\lambda$ decays to zero a little more slowly than in the typical least-squares smoothing spline problem.

\begin{theorem}
\label{Thm:SMSP2}
Assume (A1)-(A5), $\lambda \to 0 $ and $n \lambda^{3/2m -1/4m^2} \to \infty$, as $n \to \infty$. Then there exists a sequence of M-type smoothing splines $\widehat{f}_n$ satisfying
\begin{equation*}
||\widehat{f}_n - f_o||_{m, \lambda}^2 = O_P\left( \lambda + (n\lambda^{1/2m})^{-1} \right).
\end{equation*}
\end{theorem}
The limit requirements of Theorem~\ref{Thm:SMSP2} replace the least-squares requirements $\lambda \to 0$ and $n \lambda^{1/2m} \to \infty$, as $n \to \infty$. For $\lambda \asymp n^{-2m/(2m+1)}$ these conditions are met and we are lead to $||\widehat{f}_n - f_o||_{2}^2 = O_P(n^{-2m/(2m+1)})$, which is the optimal rate of convergence for $f_o \in \mathcal{W}^{m,2}([0,1])$ \citep{Stone:1982}. Thus, a broad class of smoothing spline estimators is theoretically optimal. Moreover, for such $\lambda$, the Sobolev embedding theorem \eqref{eq:SMSP6} allows us to deduce that $||\widehat{f}_n - f||_{\infty} = O_P(n^{(1-2m)/(2m+1})$, which implies that convergence can be made uniform.

Corollary \ref{Cor:SMSP1} below establishes (optimal) rates of convergence for the derivatives $\widehat{f}_n^{(j)}$ for $j = 1, \ldots, m-1$ and tightness of $\widehat{f}_n^{(m)}$ in the $\mathcal{L}^2([0,1])$ metric.

\begin{cor}
\label{Cor:SMSP1}
Assume the conditions of Theorem~\ref{Thm:SMSP2} hold. Then, for any $\lambda \asymp n^{-2m/(2m+1)}$ the M-type sequence $\widehat{f}_n$ of Theorem \ref{Thm:SMSP2} satisfies
\begin{equation*}
|| \widehat{f}^{(j)}_n - f^{(j)}||^2_2 = O_P\left( n^{-2(m-j)/(2m+1)}\right).
\end{equation*}
\end{cor}
As noted previously, with the exception of the work of \citet{Eg:2009} on the least-absolute deviations smoothing spline, we are unaware of results concerning derivatives of general M-type estimates. Corollary \ref{Cor:SMSP1} serves to remedy this deficiency.

A rather interesting feature of the least-squares smoothing spline is the possibility for a bias-reduction,  under certain boundary conditions on the derivatives of $f_o$, see, e.g., \citep{Rice:1981, Eubank:1999}. This phenomenon leads to superior convergence rates and is known as super-convergence.  As Corollary \ref{Cor:SMSP2} shows, super-convergence carries over to the general M-case.

\begin{cor}
\label{Cor:SMSP2}
Assume the conditions of Theorem~\ref{Thm:SMSP2} hold and further that $f_o \in \mathcal{W}^{2m,2}([0,1])$ and $f_o^{(s)}(0) = f_o^{(s)}(1) = 0 $ for all $ m \leq  s \leq 2m-1$. Then, there exists a sequence of M-type smoothing splines $\widehat{f}_n$ satisfying
\begin{equation*}
||\widehat{f}_n - f_o||_{m, \lambda}^2 = O_P\left( \lambda^2 + (n\lambda^{1/2m})^{-1} \right).
\end{equation*}
\end{cor}
\noindent
A consequence of this corollary is that if $\lambda \asymp n^{-2m/(4m+1)}$ then $f_o$ can be estimated with an integrated squared error decaying like $n^{-4m/(4m+1)}$ asymptotically. Of course,  had we suspected that $f_o \in \mathcal{W}^{2m, 2}([0,1])$ then an appropriate modification of the penalty would also yield the rate $n^{-4m/(4m+1)}$. In this light, as \citet[p. 259]{Eubank:1999} notes, the higher rate of convergence may be viewed as a bonus of the smoothing-spline estimator for some situations where the regression function is smoother than anticipated.

We now turn to the problem of M-type smoothing splines with an auxiliary scale estimate. We aim to extend Theorem \ref{Thm:SMSP2} and its corollaries to this case under suitable assumptions on $\rho$, $\epsilon$ and $\widehat{\sigma}$. The revised set of assumptions is as follows.

\begin{itemize}
\item[(B1)] $\rho$ is a convex function on $\mathbb{R}$ with derivative $\psi$ that exists everywhere and is Lipschitz. Further, for any $\epsilon>0$ there exists $M_{\epsilon}$ such that
\begin{equation*}
\left|\psi(tx) - \psi(sx) \right| \leq M_{\epsilon} |t-s|,
\end{equation*}
for any $ t> \epsilon, s> \epsilon$ and $-\infty<x <\infty$.
\item[(B2)] For any $\alpha>0$, $\mathbb{E}\{ \psi(\epsilon_1/ \alpha) \} = 0$ and
\begin{equation*}
\mathbb{E}\left\{\psi\left( \frac{\epsilon_1}{\alpha} + t\right) \right\} = \xi(\alpha) t + o(t), \quad \text{as} \quad t \to 0,
\end{equation*}
for some function $\xi(\alpha)$ such that $0<\inf_{|\sigma - \alpha| \leq \delta} \xi(\alpha) \leq \sup_{|\sigma - \alpha| \leq \delta} \xi(\alpha)< \infty$ for some $\delta>0$.
\item[(B3)] $n^{m/(2m+1)}(\widehat{\sigma}-\sigma) = O_P(1)$,
for some $\sigma>0$ that need not be the standard deviation of $\epsilon_1$.
\item[(B4)] (A5).
\end{itemize}

The above conditions essentially require a somewhat smoother $\rho$-function and a well-behaved scale. Condition (B1) is borrowed from \citep{He:1995} and requires the score function to change slowly in the tail. It automatically implies its boundedness. If $\psi$ is differentiable, this condition is satisfied whenever $\sup_{x}|x \psi^{\prime}(x)| < \infty$. Condition (B2) extends the linearization imposed by (A4) to a neighbourhood around $\sigma$. These conditions are not stringent and are satisfied by, e.g., the Huber score function with  $\xi(\alpha) = 2 F(k/\alpha) -1$.

Assumption (B4) imposes a minimal rate of convergence for $\widehat{\sigma}$, namely the optimal non-parametric rate of convergence. To the best of our knowledge, this is the weakest condition that is currently available in the literature, as it can be fulfilled not only by scale estimates obtained from linear combinations of the $Y_{i}$, i.e., pseudo-residuals, but also by scale estimates that are constructed  from the residuals of an initial regression estimator. The former class of scale estimates normally converges at root-n rate, see, e.g., \citep{Cunningham:1991, Ghement:2008}, while the latter class generally converges at the optimal non-parametric rate required herein. It is important to note that the above conditions do not include any moment assumptions on $\epsilon$, which are undesirable from the point of view of robustness. This needs to be contrasted with the first moment assumed by \citet{Cunningham:1991} and \citet{Oh:2007} or the second moment assumed by \citet{He:1995}.

We denote the sequence of M-estimates with auxiliary scale $\widehat{f}_{n ,\widehat{\sigma}}$. Theorem \ref{Thm:SMSP2} now generalizes to Theorem \ref{Thm:SMSP3} below.

\begin{theorem}
\label{Thm:SMSP3}
Assume (B1)-(B4) and further that $\lambda \to 0 $ and $ n \lambda^{1/m} \to \infty$, as $n \to \infty$. Then there exists a sequence of M-type smoothing splines $\widehat{f}_{n, \widehat{\sigma}}$ satisfying
\begin{equation*}
||\widehat{f}_{n, \widehat{\sigma}} - f_o||_{m, \lambda}^2 = O_P\left( \lambda + (n\lambda^{1/2m})^{-1} \right).
\end{equation*}
\end{theorem}
\noindent	
Since the proofs of Corollaries \ref{Cor:SMSP1} and \ref{Cor:SMSP2} regarding optimal estimation of derivatives and bias reduction do not depend on the existence of a scale estimate, they carry over immediately to this setting.

As requested by a reviewer, we may also seek conditions ensuring the extension of the above results to the case of a random design, that is, to the case where $t_i, \ i =1, \ldots, n$, are i.i.d. random variables with distribution $G$. The main complication in this setting is that (A5) may fail to hold due to the comparatively lower rate of convergence of the empirical distribution to its population counterpart in the uniform metric. Nevertheless, inspection of the proofs reveals that such an extension is indeed possible after the modification of (A5) given below. To better indicate the randomness of the design points we now write $T_1, \ldots, T_n$ instead of the previous $t_1, \ldots, t_n$.
\begin{itemize}
\item[(A5$^{\prime}$)] The design variables $T_1, \ldots, T_n$ are independent and identically distributed with Lebesgue density $g$ that is bounded away from zero and infinity on $[0,1]$, and are also independent of the errors $\epsilon_1, \ldots, \epsilon_n$.
\end{itemize}
The assumption of a bounded density goes back at least to \citep{Stone:1985} and has been extensively used in non-parametric regression. With this assumption in place, Theorem~\ref{Thm:SMSP4} below now replaces Theorem~\ref{Thm:SMSP2}.

\begin{theorem}
\label{Thm:SMSP4}
Assume (A1)-(A4), (A5$^{\prime}$), $\lambda \to 0 $ and $n \lambda^{3/2m - 1/4m^2} \to \infty$, as $n \to \infty$. Then there exists a sequence of M-type smoothing splines $\widehat{f}_n$ satisfying
\begin{equation*}
||\widehat{f}_n - f_o||_{m, \lambda}^2 = O_P\left( \lambda + (n\lambda^{1/2m})^{-1} \right).
\end{equation*}
\end{theorem}

An extension of Theorem~\ref{Thm:SMSP3} to the present setting is now completely analogous under the understanding that (B5) there is replaced with condition (A5$^{\prime}$), as are the extensions of Corollary 1 and Corollary 2. We omit the details.

\section{Computation and smoothing parameter selection}
\label{sec:SMSP3}

As discussed in Section~\ref{sec:SMSP2}, there exists at least one solution of \eqref{eq:SMSP3} in the space of natural splines of order $2m$ with knots at the unique $t_{i}$.  Thus we may restrict attention to the linear subspace of natural splines for the computation of the estimator. Here, we assume that the scale is known and equal to one; if that is not true, the scale can be absorbed into the $\rho$-function and all the arguments of this section still apply. Assume for simplicity that all $t_{i}$ are distinct and let $a = \min_{i} t_{i}>0$ and $b= \max_{i} t_{i}<1$. Then the natural spline $f$ has $n$ interior knots and we may write
\begin{equation}
\label{eq:SMSP9}
f(t) = \sum_{k = 1}^{n + 2m} f_k B_{k}(t),
\end{equation}
where $f_k$ are scalar coefficients and the $B_k(\cdot)$ are the B-spline basis functions of order $2m$ supported by the knots at $t_{i}$.  It is well-known that the space of natural splines with $n$ interior knots is $n$-dimensional on account of the boundary conditions that force the natural spline to be a polynomial of order $m$ outside $[a,b]$, but for the moment we operate in the larger $(n + 2m)$-dimensional spline subspace. This is computationally convenient due to the local support and numerical stability of B-splines. As we explain below, the penalty automatically imposes the boundary conditions \citep[see also][pp. 161-162]{Hastie:2009}.

Working in the unrestricted spline subspace, the solution to \eqref{eq:SMSP3} or \eqref{eq:SMSP4} may be written as $\widehat{f}_n = \sum_{k = 1}^{n + 2m} \widehat{f}_k B_{k}(t)$ where $\boldsymbol{\widehat{f}} = (\widehat{f}_1, \ldots, \widehat{f}_{n+2m})^{\top}$ satisfies
\begin{align}
\label{eq:SMSP10}
\boldsymbol{\widehat{f}} = \argmin_{\boldsymbol{f} \in \mathbb{R}^{n+2m}} \left[ \frac{1}{n} \sum_{i=1}^n \rho\left(Y_i - \mathbf{B}_i^{\top} \mathbf{f}\right) + \lambda \mathbf{f}^{\top} \mathbf{P} \mathbf{f} \right],
\end{align}
with  $\mathbf{B}_{i} = (B_1 (t_{i}), \ldots, B_{n+2m}(t_{i}) )^{\top}$ and $\mathbf{P} = \langle B_k^{(m)}, B_l^{(m)} \rangle_2,  k,l = 1, \ldots, n+2m$. Initially, it may seem that this formulation ignores the boundary constraints that govern natural splines but it turns out the penalty term automatically imposes them. The reasoning is as follows: if that were not the case, it would always be possible to find a $2m$th order natural interpolating spline of the form of \eqref{eq:SMSP9} that leaves the first term in \eqref{eq:SMSP3} unchanged, but due to it being a polynomial of order $m$ outside of $[a,b] \subset [0,1]$ the penalty semi-norm  would be strictly smaller. Hence, any minimizer in the form of \eqref{eq:SMSP9} incorporates the boundary conditions.

The solution to \eqref{eq:SMSP10} may be expediently found through minor modification of the penalized iteratively reweighted least-squares algorithm found in, e.g., \citep{Maronna:2011}. The algorithm consists of solving a weighted penalized least-squares problem at each iteration until convergence, which is guaranteed irrespective of the starting values and yields a stationary point of \eqref{eq:SMSP3}, under mild conditions on $\rho$ that include the boundedness of $\rho^{\prime}(x)/x$ near zero~\citep{Huber:2009}. For the quantile loss, this condition fails on account of the kink at the origin. Nevertheless, the easily implementable recipe of \citep{Nychka:1995} may be used in order to obtain an approximate solution of \eqref{eq:SMSP3}. In particular, in the algorithm the loss function can be replaced by the smooth approximation
\begin{equation*}
\widetilde{\rho}_{\alpha}(x) = \begin{cases} \rho_{\alpha}(x) &  |x| \geq \epsilon \\
\alpha x^2/\epsilon & 0 \leq x < \epsilon \\
(1-\alpha)x^2/\epsilon & -\epsilon <  x \leq 0,
\end{cases}
\end{equation*}
for some small $\epsilon>0$. The objective function $\widetilde{\rho}_{\alpha}$  allows to easily calculate an approximate quantile smoothing spline estimate with this algorithm, without the need of a more computationally burdensome quadratic program. Whenever such a modification of the objective function is not feasible, we recommend utilizing a convex-optimization program in order to identify a minimizer of \eqref{eq:SMSP10}, as given, for example, by \citep{Fu:2020}.

In order to select the smoothing parameter parameter $\lambda$ we propose to use the weighted generalized cross-validation criterion proposed by \citep{Cunningham:1991}. That is, we propose to select $\lambda$ as the minimizer of 
\begin{equation*}
\GCV(\lambda) = \frac{n^{-1} \sum_{i=1}^n W_i(\widehat{f}_n) |r_i(\widehat{f}_n)|^2}{|1-n^{-1}\Tr \mathbf{H}(\lambda)|^2},
\end{equation*}
where $r_i(\widehat{f}_n) = Y_i - \widehat{f}_n(t_i)$, $W_i(\widehat{f}_n) = \psi(r_i(\widehat{f}_n))/r_i(\widehat{f}_n)$, $i=1, \ldots, n,$ and $\mathbf{H}(\lambda)$ is the pseudo-influence matrix obtained upon convergence. Throughout the simulation experiments and real-data examples to follow we have adopted a two-step approach in order to identify the minimizer of $\GCV(\lambda)$. First, we have determined the approximate location of the minimizer by evaluating $\GCV(\lambda)$ on a grid and afterwards employed a numerical optimizer based on golden section search and parabolic interpolation \citep{Nocedal:2006}. Such a hybrid approach is often advisable due to the possible local minima and near-flat regions of the GCV criterion, particularly for non-smooth loss functions.

\section{A Monte-Carlo study}
\label{sec:SMSP4}

In our simulation experiments we compare the performance of the Huber-type smoothing spline with tuning parameter equal to $1.345$ corresponding to 95\% efficiency in the location model, the least-absolute deviation smoothing spline and the least-squares smoothing spline in a variety of shapes of the regression function and tails of the error. For the Huber-type estimator we compute the scale in two ways, reflecting the asymptotic considerations of Section~\ref{sec:SMSP2}. 

On the one hand, we use a robust scale constructed from consecutive differences of the response variables, as proposed by \citet{Ghement:2008}. In particular, we use the M-scale estimator $\widehat{\sigma}$ obtained as the solution of
\begin{align*}
\frac{1}{n-1} \sum_{i=1}^{n-1} \rho \left( \frac{Y_{i+1}- Y_i}{2^{1/2} \widehat{\sigma}} \right) = \frac{3}{4},
\end{align*}
where $\rho$ is the Tukey loss function with tuning parameter equal to $0.704$. The constants $2^{1/2}$ and $3/4$ ensure Fisher-consistency of the estimator at the Gaussian distribution and maximal breakdown value, respectively. On the other hand,  call $\mathbf{r}$ the residual vector that results from an initial $L_1$-type smoothing spline fit to the data. A robust scale that converges slower than root-$n$ is
\begin{equation*}
\widehat{\sigma}_{PR} = \tau(\mathbf{r}),
\end{equation*}
where $\tau(\cdot)$ stands for the $\tau$-scale introduced by \citet{yohai1988high}. In our implementation we used the Tukey bisquare loss function and cut-off constants equal to $4.5$ and $3$ for the biweighting of the mean and the $\rho$-function respectively, as recommended by \citet{Maronna:2002}. As discussed in the introduction, the least-absolute deviations smoothing spline does not require an auxiliary scale estimate. We denote the resulting smoothing splines by $\widehat{f}_{HPS}$ and $\widehat{f}_{HPR}$ respectively, while we use $\widehat{f}_{LAD}$ and $\widehat{f}_{LS}$ as abbreviations for the least-absolute deviations smoothing spline and the least-squares smoothing spline respectively.

We investigate the performance of the estimators in the regression model $Y_i = f(t_i) + \epsilon_i$ where $t_i = (i-1/2)/n$ and $f_o$ is each of the following three functions
\begin{enumerate}
\item $f_1(t) = \cos(2 \pi t)$,
\item $f_2(t) = 1/(1+\exp(-20(t-0.5)))$,
\item $f_3(t) = \sin(2 \pi t) + e^{-3(t-0.5)^2}$.
\end{enumerate}
All three regression functions are smooth; the first regression function is bowl-shaped, the second is a sigmoid and the third is essentially a shifted sinusoid with more variable slope. We shall estimate these functions with cubic smoothing splines, that is, with $m=2$. In order to assess the effect of outliers on the estimates different distributions for the error term were considered. Other than the standard Gaussian distribution, we have complemented our set-up with a t-distribution with $3$ degrees of freedom, a right-skewed t-distribution with 3 degrees of freedom and non-centrality parameter equal to $2$, a mixture of mean-zero Gaussians with standard deviations equal to $1$ and $9$ and weights equal to $0.85$ and $0.15$ respectively, as well as Tukey's Slash distribution defined as the quotient of independent standard Gaussian and uniform random variables. 

All computations were carried out in the freeware \textsf{R} \citep{R}. For ease of comparison all smoothing spline estimators were computed using custom-made functions implementing the method outlined in Section~\ref{sec:SMSP3}. For all smoothing spline estimators the penalty parameter was selected via the GCV criterion given in Section~\ref{sec:SMSP3}. A link to the implementations used in this experiment is given at the end of this paper. The mean-squared-errors of the experiment are summarized in Table \ref{tab:SMSP1} for sample sizes of 60 and 1000 replications.
\begin{table}[h]
\centering
\begin{tabular}{lcccccccccc}
& & \multicolumn{2}{c}{$\widehat{f}_{HPS}$} & \multicolumn{2}{c}{$\widehat{f}_{HPR}$} & \multicolumn{2}{c}{$\widehat{f}_{LAD}$} & \multicolumn{2}{c}{$\widehat{f}_{LS}$} \\ \\[-2ex]
$f$ & Dist. & Mean & SE & Mean & SE & Mean & SE & Mean & SE  \\ \\
\multirow{4}{*}{$f_1$} & Gaussian & 0.087  & 0.002 & 0.086 & 0.002 & 0.114 & 0.002 & 0.082 & 0.002 \\
 & $t_3$ & 0.125 &  0.003 & 0.125 & 0.003 & 0.137 & 0.003 & 0.216 & 0.016 \\
 & $st_{3,2}$ & 0.771 &  0.008 & 0.794 & 0.008 & 0.649 & 0.008 & 1.108 & 0.014 \\
 & M. Gaussian &  0.138 & 0.004 & 0.139 & 0.003 & 0.146 & 0.004 & 0.719 & 0.020 \\ 
 & Slash & 0.520 & 0.013 & 0.515 & 0.012 & 0.389 & 0.009 & 766.6 & 352.9 \\ \\
\multirow{4}{*}{$f_2$} & Gaussian &  0.064 & 0.001 & 0.064 & 0.001 & 0.081 & 0.002  & 0.065 & 0.001 \\	
 & $t_3$ & 0.091  & 0.002 & 0.090 &  0.002 & 0.101 & 0.002 & 0.153 & 0.056 \\
 & $st_{3,2}$ & 0.767  & 0.008 & 0.765 &  0.008 & 0.635 & 0.007 & 1.107 & 0.014 \\
 & M. Gaussian & 0.097  & 0.004  & 0.095 & 0.004 & 0.106 & 0.003 & 0.562 & 0.019	 \\ 
 & Slash & 0.345 & 0.011 & 0.331  & 0.011 & 0.286 & 0.008 & 7997 &  6934 \\ \\
\multirow{4}{*}{$f_3$} & Gaussian &  0.117  & 0.003  & 0.117 & 0.003 & 0.166 & 0.002 & 0.099 &  0.002 \\
 & $t_3$ &  0.182 &  0.003 & 0.183 & 0.003  & 0.192 & 0.003 & 0.234 & 0.005 \\
 & $st_{3,2}$ &  0.793 &  0.008 & 0.816 & 0.008  & 0.721 & 0.008 & 1.134 & 0.015 \\
 & M. Gaussian & 0.203 & 0.004 & 0.204 & 0.004 & 0.203 & 0.003 & 0.734 & 0.019 \\ 
 & Slash & 0.518  & 0.011 &  0.512 & 0.010 & 0.418 & 0.008 &  1e+08 & 1e+08	 \\ 
\end{tabular}
\caption{Means and standard errors of 1000 MSEs with $n=60$ of the Huber-type estimator with preliminary scale, the Huber-type estimator with scale computed from robust regression residuals, the least-absolute deviations estimator and the least-squares estimator.}
\label{tab:SMSP1}
\end{table}
\normalsize

The results in Table \ref{tab:SMSP1} indicate the extreme sensitivity of the least-squares estimator to even mild deviations from the ideal assumptions. In particular, the  regular t-distribution with 3 degrees of freedom and the mixture-Gaussian distribution have finite second moments and first moments equal to zero, so that the least-squares assumptions are technically fulfilled. Nevertheless, in all three of our examples the performance of the least-squares estimator markedly deteriorates as one moves away from the Gaussian ideal. By contrast, the robust estimators almost match the performance of the least-squares estimator in ideal data and exhibit a large degree of resistance in the t- and mixture Gaussian distributions. Quite notably, the robust estimators seem somewhat vulnerable to asymmetric contamination, although they still outperform the least-squares estimator in this setting. Their performance also deteriorates with the heavy-tailed Slash distribution but clearly not to the same extent as the performance of the least-squares estimator, which appears to suffer a catastrophic breakdown. 

Comparing the robust estimators $\widehat{f}_{HPS}$, $\widehat{f}_{HPR}$ and $\widehat{f}_{LAD}$ in detail, it is seen that the Huber estimators outperform the least-absolute deviation estimator in the case of Gaussian, $t_3$ and mixture Gaussian errors but get outperformed in turn under skewed-t and Slash errors. These facts indicate the resistance of the least-absolute deviations estimator both with respect to asymmetric contamination and gross outliers. The Huber estimators  $\widehat{f}_{HPS}$ and $\widehat{f}_{HPR}$ overall exhibit similar performance, except in the case of heavy contamination where the latter estimator has a slight edge on account of its more robust scale. However, this advantage comes at the price of additional computational effort, as  $\widehat{f}_{HPR}$ depends on the computation of two robust smoothing spline estimators. Overall, the present experiment indicates that the computationally simple $\widehat{f}_{HPS}$ presents a viable alternative to $\widehat{f}_{LS}$ in clean data and to $\widehat{f}_{HPR}$ in mildly contaminated data.

\section{Real data example: urban air pollution in Italy}
\label{sec:SMSP5}

\begin{figure}[h]
\centering
\subfloat{\includegraphics[width = 0.495\textwidth]{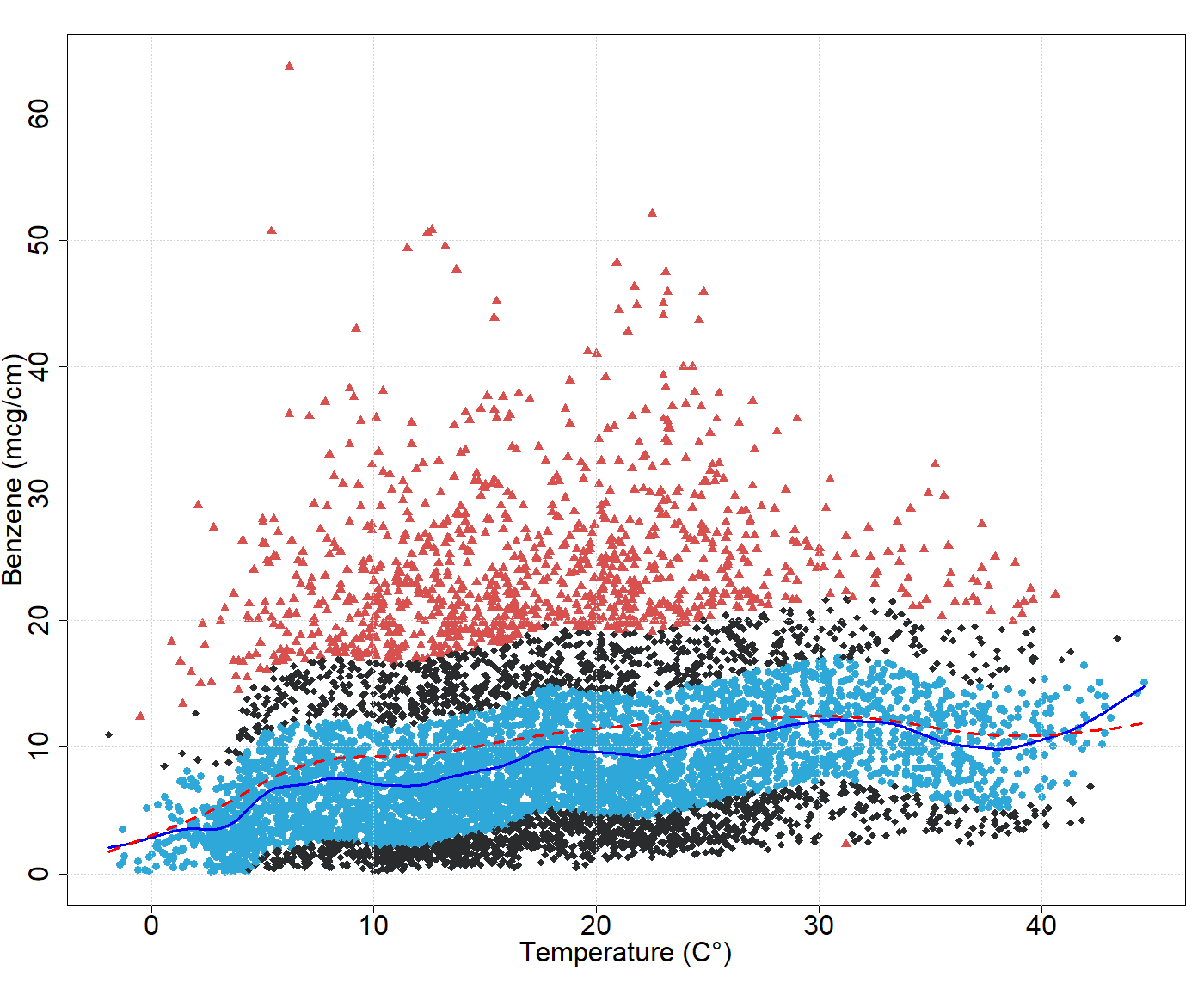}} \
\subfloat{\includegraphics[width = 0.495\textwidth]{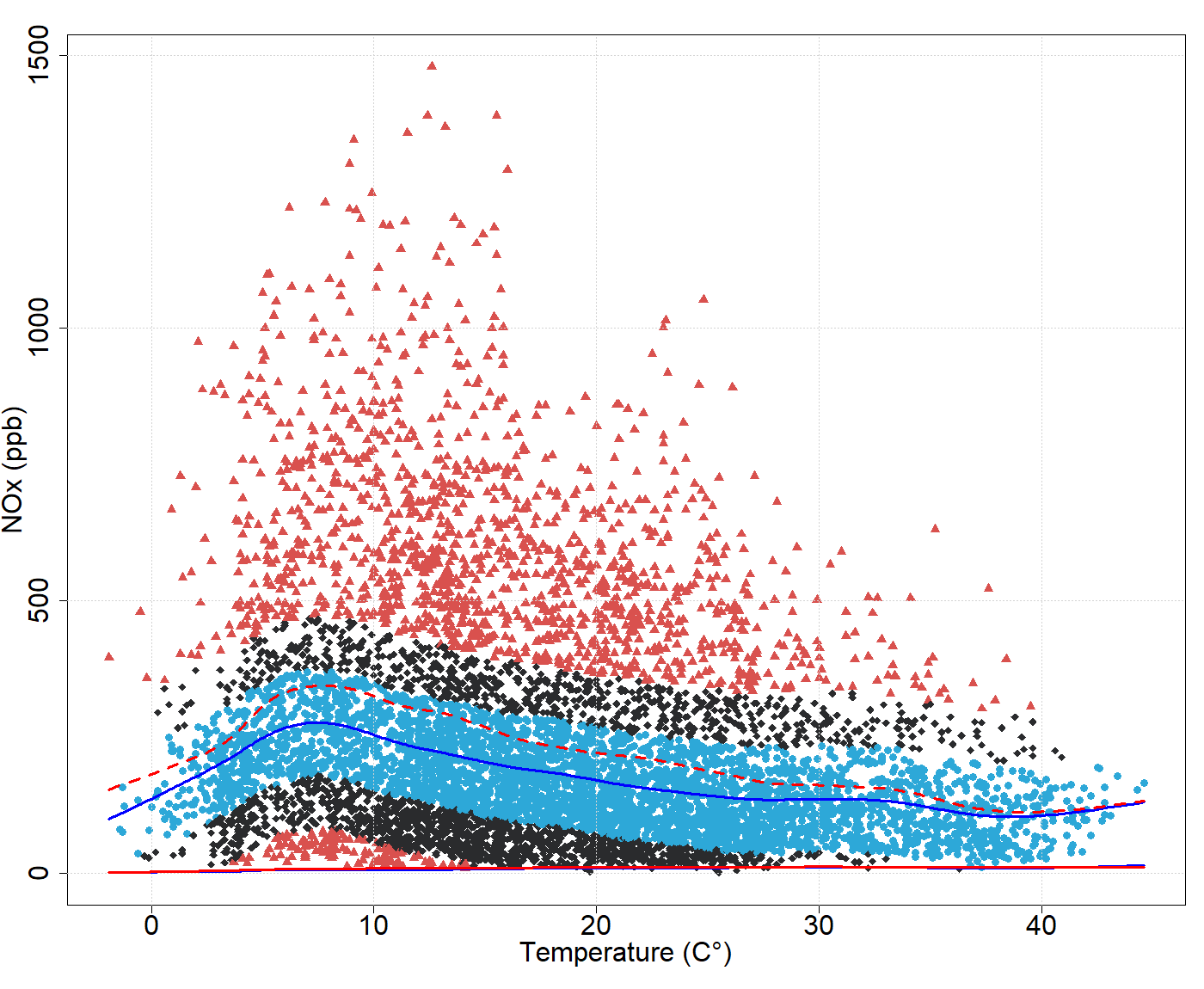}}
\caption{Left: scatter plot of benzene versus temperature. Right: scatter plot of nitrogen oxide versus temperature. The lines (\full, \dashed) correspond to the Huber and least-squares estimators respectively, while the symbols (\textcolor{BB}{$\bullet$}  $\scriptstyle{\textcolor{TS}{\blacklozenge}}$ \textcolor{DR}{$\blacktriangle$}) denote Huber weights in $(0.66, 1] ,\ (0.33, 0.66]\ \text{and}\ (0, 0.33] $ respectively. }
\label{fig:SMSP1}
\end{figure}

The present dataset consists of air pollution measurements from a gas multi-sensor device that was deployed in an undisclosed Italian city \citep{De Vito:2008}. In particular, the dataset contains 7396 instances of hourly averaged responses from an array of 5 metal oxide chemical sensors embedded in the device. In this study we will focus on the effect of temperature on the concentration of Benzene and Nitrogen oxides in the atmosphere, the full dataset being available at \url{http://archive.ics.uci.edu/ml/datasets/Air+quality}. These gases are known carcinogens and are responsible for a series of acute respiratory conditions in high concentrations, so that the ability to forecast high concentrations may be helpful. Figure~\ref{fig:SMSP1} presents the scatter plots of Benzene concentrations versus temperature and Nitrogen oxides concentration versus temperature on the left and right panel respectively.

The scatter plots indicate the presence of several atypical observations in the form of abnormally high concentrations at certain temperature ranges. Specifically, the concentration of Benzene seems more volatile in medium temperatures, while the concentration of Nitrogen oxides seems more unpredictable in lower temperatures. Computing the least-squares and Huber-type with preliminary scale smoothing splines yields the solid red and solid blue curves respectively. It may be immediately seen that although the resulting smoothing spline estimates do not qualitatively differ as to their main features, the least-squares estimate is more drawn towards these atypically high concentrations resulting in overestimation of the concentrations, particularly in the case of Nitrogen oxides. 

For a better understanding of the discrepancy between these estimates, figure~\ref{fig:SMSP1} also includes a color and shape coding of the weights $\psi(r_i)/r_i$ produced by the Huber estimator. Here, $r_i$ stands for the ith residual, i.e., $r_i = y_i - \widehat{f}(t_i)$, $i=1, \ldots, n$. These weights demonstrate the usefulness of M-estimators with bounded score functions: while the least-squares estimator assigns equal weight to all observations and is thus unduly influenced by outliers, robust M-estimators significantly downweight observations sufficiently far from the center of the data resulting in resistant estimates. These atypical observations produce large residuals, which in turn allow for their detection. We may think of numerous examples in economics, medicine and other fields where resistant estimation combined with the possibility of outlier detection can be similarly useful.

\begin{figure}[h]
\centering
\subfloat{\includegraphics[width = 0.495\textwidth]{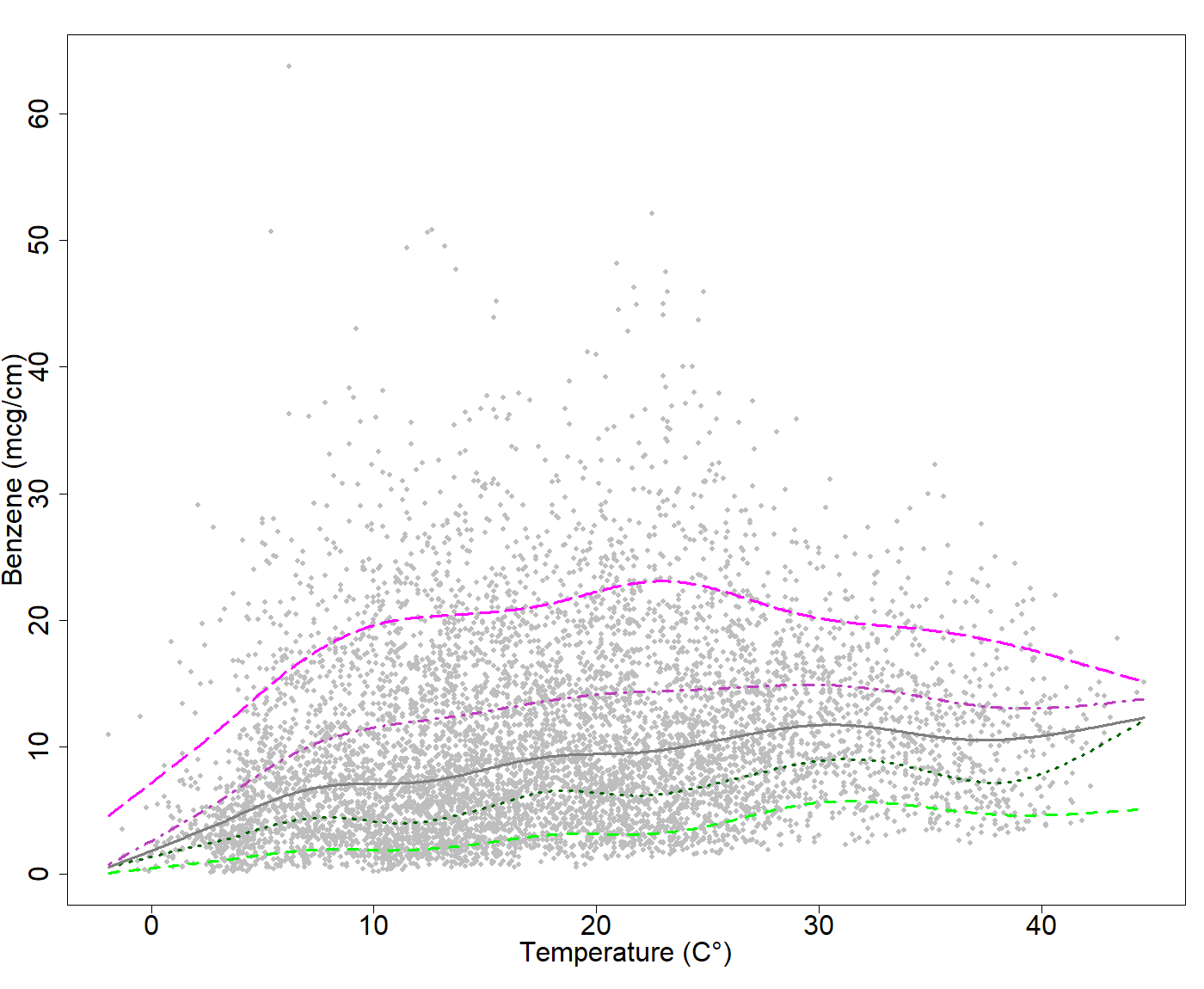}} \
\subfloat{\includegraphics[width = 0.495\textwidth]{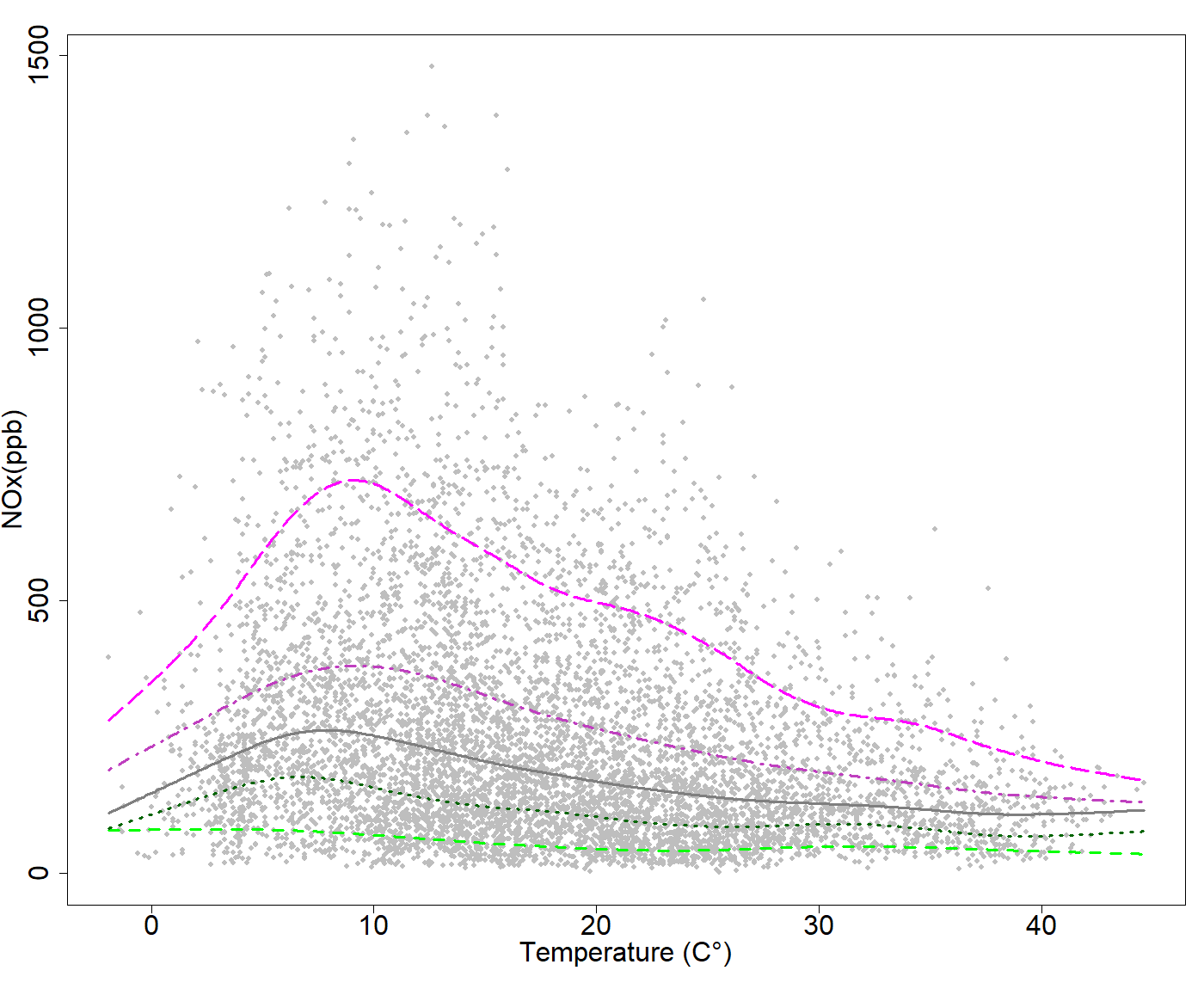}}
\caption{Left: scatter plot of benzene versus temperature. Right: scatter plot of nitrogen oxide versus temperature. The lines (\dashedd, \dotted, \fulll,  \dotdash, \denselydashed) correspond to estimated $(0.1, \ 0.3, \ 0.5, \ 0.7, \ 0.9)$-quantiles respectively.}
\label{fig:SMSP2}
\end{figure}

Lastly, to examine the characteristics of the conditional distributions of the gases given the temperature, one may wish to obtain nonparametric estimates of the conditional quantiles. Figure \ref{fig:SMSP2} presents quantile smoothing spline estimates for the 10th, 30th, 50th, 70th and 90th percentiles. The estimates serve to confirm the extreme heavy-tailedness of the conditional distributions for medium and medium-low temperatures, casting doubt on the suitability of the least-squares estimator.

\section{Concluding remarks}
\label{sec:SMSP6}

The asymptotic results of this paper indicate that there is little theoretical difference between least-squares smoothing splines and general M-type smoothing splines derived from convex but possibly non-smooth objective functions. In particular, under general conditions, M-type smoothing splines enjoy the same rates of convergence for the regression function as well as its derivatives. Furthermore, the presence of reasonable auxiliary scale estimates does not diminish these rates. In practice, M-type smoothing splines may be efficiently computed with the convenient B-spline representation and well-established iterative algorithms and may be used to good effect in either regular or contaminated datasets. Thus, this broad family of estimators provides a valuable tool for the applied scientist.

We believe that at least two extensions of the present paper would be of great interest to theoreticians and practitioners alike. The first concerns the problem of estimating the smoothing parameter $\lambda$, which was barely touched upon in this paper. We conjecture that the robust generalized cross-validation criterion proposed by \citet{Cunningham:1991} would yield the optimal rate of decay for $\lambda$ but formal verification is required. Another useful extension would be to the case of dependent errors, which arise, for example, in mean-estimation of discretely sampled functional data. We firmly believe that M-type smoothing spline estimators with repeated measurements would still enjoy the optimal rates, as derived by \citet{Cai:2011}, while also providing a considerably safer estimation method in the presence of atypical observations.

\section*{Acknowledgements}

The present revised version contains a correction in the statement and proof of Theorem~\ref{tab:SMSP1}. The author is grateful to two anonymous referees and an associate editor for valuable comments and suggestions that greatly improved the paper both with respect to content and accessibility. He is also indebted to Stefan Van Aelst for helpful conversations. Support from grant C16/15/068 of KU Leuven is gratefully acknowledged.

\section*{Software availability}

\textsf{R}-scripts reproducing the simulation experiments and real-data examples are available in \url{https://github.com/ioanniskalogridis/Smoothing-splines}.

\bibliographystyle{spbasic}      


\section*{Appendix: Proof of Theorem~\ref{Thm:SMSP1}}

Let $L_n(f)$ denote the objective function, that is,
\begin{align*}
L_n(f) = \frac{1}{n} \sum_{i=1}^n \rho(Y_i-f(t_i))+\lambda \|f^{(m)}\|^2_2.
\end{align*}
Let us first show that $L_n: \mathcal{W}^{m,2}([0.1]) \to \mathbbm{R}$ is weakly lower semicontinuous. Let $\{f_k\}_k$ denote a sequence of functions in $\mathcal{W}^{m,2}([0,1])$ converging weakly to some $f\in \mathcal{W}^{m,2}([0,1])$. By the reproducing property,
\begin{align*}
\lim_{k\to \infty} f_k(t_i) = \lim_{k \to \infty} \langle f_k, \mathcal{R}_{m,\lambda}(t_i, \cdot) \rangle_{m,\lambda} = \langle f, \mathcal{R}_{m,\lambda}(t_i, \cdot) \rangle_{m, \lambda} = f(t_i),
\end{align*}
for all $t_i$. On the other hand, the functional $f \to \|f^{(m)}\|^2_2$ is convex and lower semicontinuous and therefore also weakly lower semicontinuous. Putting these facts together, the lower semicontinuity of $\rho$ yields 
\begin{align*}
L_n(f) & = \frac{1}{n} \sum_{i=1}^n \rho(Y_i-f(t_i)) + \lambda \|f^{(m)}\|^2_2
\\ & \leq \frac{1}{n} \sum_{i=1}^n \liminf_{k \to \infty} \rho(Y_i-f_k(t_i)) + \lambda \liminf_{k \to \infty} \|f_k^{(m)}\|^2_2
\\ & \leq  \liminf_{k \to \infty} \frac{1}{n} \sum_{i=1}^n \rho(Y_i-f_k(t_i)) + \lambda \liminf_{k \to \infty} \|f_k^{(m)}\|^2_2
\\ & \leq \liminf_{k \to \infty} L_n(f_k),
\end{align*}
establishing the weak lower semicontinuity of the objective function as a whole.

Now note that $L_n(f)$ is bounded from below by $0$ for all $f \in \mathcal{W}^{m,2}([0,1])$, hence its infimum over $\mathcal{W}^{m,2}([0,1])$ is finite. Let $\{f_k\}_k$ denote a minimizing sequence, that is,
\begin{align*}
\lim_{k \to \infty} L_n(f_k) = \inf_{f \in \mathcal{W}^{m,2}([0,1])}L_n(f).
\end{align*}
By Taylor's theorem with integral remainder we may write
\begin{align}
\label{eq:A1}
f_k(t) = p_k(t) +  \int_{0}^{1} \frac{(t-x)_{+}^{m-1}}{(m-1)!} f_k^{(m)}(x) dx,
\end{align}
where $p_k$ is a polynomial of order $m$ and $(x)_{+}^{m-1} = x^{m-1} \mathcal{I}(x \geq 0)$.

Let $T: \mathcal{L}^2([0,1]) \to \mathcal{C}{([0,1])}$ denote the integral operator 
$g(t) \mapsto \int_{0}^t (t-x)^{m-1}/(m-1)! g(x) dx$. We claim that $T$ is compact. Indeed let $\{f_n\}_n$ denote a bounded sequence in $\mathcal{L}^2([0,1])$ with $\sup_{n} \|f_n\|_2 \leq M$, say. By the Schwarz inequality 
\begin{align*}
\sup_{n} |Tf_n(t)| \leq \frac{|t|^{(2m-1)/2}}{(2m-1)^{1/2}(m-1)!} \sup_n \|f_n\|_2 \leq \frac{M}{(2m-1)^{1/2} (m-1)!}.
\end{align*}
Further, for $m \geq 2$ by uniform continuity for every $\epsilon>0$ there exists a $\delta>0$ such that $|(t-x)_{+}^{m-1}-(t^{\prime}-x)_{+}^{m-1}| < \epsilon$ uniformly in $x \in [0,1]$ whenever $|t-t^{\prime}| < \delta$. It follows that 
\begin{align*}
|Tf_n(t) - Tf_n(t^{\prime})| \leq \frac{\epsilon}{(m-1)!}\int_{0}^1 |f_n(x)| dx \leq \frac{\epsilon}{(m-1)!}\|f_n\|_2 \leq \frac{\epsilon}{(m-1)!}M.
\end{align*}
For $m=1$, $Tf_n(t) = \int_{0}^t f_n(x) dx$, which is also equicontinuous for $\|f_n\| \leq M$. Hence for bounded sequences $\{f_n\}_n$, $\{Tf_n\}_n$ is bounded and equicontinuous. By the Arzel\`a-Ascoli theorem there exists a subsequence $\{T f_{n_k}\}_k$ that converges uniformly to some element of $\mathcal{C}([0,1])$. Conclude that $\{Tf_n\}_n$ is (relatively) compact for bounded $\{f_n\}_n$. This establishes the compactness of $T$.

Now assume without loss of generality $L_n(f_k) \leq L_n(f_1)$ for all $k \geq 1$, passing to a subsequence if necessary. It follows that $\|f^{(m)}\|^2_2 \leq \lambda^{-1} L_n(f_1)$ so that $\|f^{(m)}\|^2_2$ is bounded. By the reflexivity of $\mathcal{L}^2([0,1])$ we may assume that $f_k^{(m)}$ converges weakly to some function $\phi$, again passing to a subsequence if necessary. By the weak lower semicontinuity of the norm we have 
\begin{align*}
\|\phi\|^2_2 \leq \liminf_{k \to \infty} \|f_{k}^{(m)}\|_2^2<\infty,
\end{align*}
as weakly convergent sequences are bounded. Since $T$ is compact we also have $\lim_k \sup_{x} |Tf_k^{(m)}(x)-T\phi(x)| = 0$. On the other hand, we also have
\begin{align*}
\sum_{i=1}^n \rho(Y_i - f_k(t_i)) \leq n L_n(f_1),
\end{align*}
so that the $f_k(t_i)$ are bounded. By the Bolzano-Weierstrass theorem, we may assume that 
\begin{align*}
(f_k(t_1), \ldots, f_k(t_n)),
\end{align*}
converges in $\mathbb{R}^n$ to some vector $v$, passing to a subsequence if necessary. From \eqref{eq:A1} we may now deduce that
\begin{align*}
\lim_{k \to \infty} p_k(t_i) = v_i - T \phi(t_i), \quad i=1, \ldots, n,
\end{align*}
that is, the limit exists for each $t_i$.  Let $P$ denote the space of polynomials of order $m$ on $[0,1]$ and note that $P$ is finite-dimensional (for example, with the basis $1, t, \ldots, t^{m-1}$). Let $\ell_i(p) := p(t_i)$ denote the evaluation functionals and on $P$ define the norm
\begin{align*}
S(p) := \sum_{i=1}^n |\ell_i(p)|.
\end{align*}
It is easy to check that $S$ is a norm. It is homogeneous, satisfies the triangle inequality and $S(p) = 0$ implies that $p \in P$ has at least $n$ roots. Since, by assumption, $n \geq m$ it must be that $p$ is the zero polynomial. Now note that, since the limit at each $t_i$ exists, $\{p_k\}_k$ is Cauchy under $S(p)$ and as $P$ is finite-dimensional, hence complete, $\{p_k\}_k$ converges to a unique polynomial $p$.

Defining $\psi := p + T \phi$ and noting that $\psi \in \mathcal{W}^{m,2}([0,1])$ and $(T\phi(t))^{(m)} = \phi(t)$ we obtain
\begin{align*}
L_n(\psi) \leq \liminf_{k \to \infty} L_n(f_k) = \inf_{f \in \mathcal{W}^{m,2}([0,1])} L_n(f).
\end{align*}
This shows that the infimum of $L_n(f)$ is attained in $\mathcal{W}^{m,2}([0,1])$. The proof is complete.

\end{document}